\documentclass[12pt,reqno]{amsart}
\usepackage[italian,english]{babel}

 \usepackage{geometry}
\usepackage[latin1]{inputenc}
\usepackage[italian,english]{babel}
\usepackage{amsmath, amsfonts, amsthm, xcolor}
\usepackage{paralist}
\numberwithin{equation}{section}
\usepackage{mathtools, mathabx}

\newtheorem{thm}{Theorem}[section]

\theoremstyle{definition}
\newtheorem{rem}[thm]{Remark}
\theoremstyle{remark}

\newcommand{\R}{\mathbb{R}}

\newcommand{\de}{\partial}

\newcommand{\haus}{\mathcal{H}^{n-1}}

\title[Two inequalities for the  first Robin  eigenvalue of the Finsler Laplacian.]{
Two inequalities for the  first Robin  eigenvalue of the Finsler Laplacian.
}
\author{
Giuseppina di Blasio}
\address{Universit\`a degli Studi della Campania ``Luigi Vanvitelli'', viale Lincoln, 5 - 81100 Caserta, Italia. }\email[G. di Blasio]{
giuseppina.diblasio@unicampania.it}
\author{
Nunzia Gavitone} %
\address{Universit\`a degli Studi di Napoli Federico II, Dipartimento di Matematica e Applicazioni ``R. Caccioppoli'', Via Cintia, Monte S. Angelo - 80126 Napoli, Italia.}
\email[N. Gavitone]{nunzia.gavitone@unina.it}

\begin{document}

\maketitle


\begin{abstract} Let $\Omega \subset \R^n$, $n \ge2,$ be a bounded connected, open set with Lip\-schitz boundary. Let $F$ be a suitable norm in $\R^{n}$ and let  $\Delta_{F}u=\text{div} \left(F_{\xi}(\nabla u)F(\nabla u)\right)$ be the so-colled Finsler Laplacian, with $u \in H^{1}(\Omega)$. In this paper we prove two inequalities for $\lambda_{F}(\beta, \Omega)$, the first eigenvalue of $\Delta_{F}$ with Robin boundary conditions  involving a   positive function $\beta(x)$. As a consequence of our result we obtain the asymptotic  behavior of $\lambda_{F}(\beta, \Omega)$ when  $\beta$ is a positive constant which goes to zero.
\medskip


\end{abstract}

\section{Introduction}
Let $\Omega \subset \R^n$, $n \ge 2$, be a bounded, connected, open set with Lipschitz boundary.\\
Let  $F\colon \R^{n}\mapsto  [0,+\infty[$, be a $C^2(\R^n\setminus \{0\})$, convex  and  positively 1-homogeneous function such that
\begin{equation}
\label{eq:lin}
a|\xi| \le F(\xi)\le b |\xi|,\quad \xi \in  \R^{n},
\end{equation}
for some positive constants $a$ and $b$. Throughout the paper we will assume that $F(\xi)$ is strongly convex, that is  \begin{equation}
\label{strong}
[F^{2}]_{\xi\xi}(\xi)\text{ is positive definite in } \R^{n}\setminus\{0\}.
\end{equation}
In what follows we  assume that  $\beta\colon \partial\Omega \to ]0,+\infty[$ is a continuous  function and we define
\begin{equation}\label{m}
m:=\displaystyle\int_{\partial \Omega}\beta(x) F(\nu) \,d \mathcal H^{n-1}>0,
\end{equation}
where $\nu$ is the unit outer normal to the boundary and $d \mathcal H^{n-1}$ denotes the $(n-1)$- dimensional Hausdorff measure. Let us consider the following Robin eigenvalue problem
\begin{equation}
\label{prob}
\left\{
\begin{array}{ll}
-\Delta_Fu =\lambda_F(\beta,\Omega) u & \text{in }\Omega \\
F(\nabla u)F_{\xi}(\nabla u)\cdot \nu +\beta(x) uF(\nu)=0 & \text{on }\partial \Omega ,
\end{array}
\right.
\end{equation}
where  $u\in H^1(\Omega)$ and $$\Delta_Fu=\text{div} \left(F_{\xi}(\nabla u)F(\nabla u)\right)$$ is the so-called Finsler Laplacian.
When $F=\mathcal E$ is the Euclidean norm, $\Delta_{F}$ reduces to the classic Laplace operator. Nevertheless, it is in general a nonlinear operator and it has been studied in several papers (see for instance \cite{aflt}, \cite{bfk}, \cite{fk}, \cite{dpdfg}, \cite{dgn}).

In \cite{gt} (see also \cite{pota} for the case $\beta$ equals to a positive constant) it is proved that the first  eigenvalue of \eqref{prob} is positive, simple and has the following variational characterization
 \begin{equation}\label{lambda}
\lambda_F(\beta,\Omega)= \min_{\substack{v \in H^{1}(\Omega)\setminus \{0\}}}\displaystyle
\frac{\int_{\Omega}F^2(\nabla v)dx+\int_{\partial \Omega}\beta(x)v^2 F(\nu)d \mathcal H^{n-1}}{\int_{\Omega}v^2dx}.
\end{equation}
On the other hand,   $\lambda_F(\beta,\Omega)$ verifies a Faber-Krahn type inequality for suitable  functions $\beta(x)$.   Finally, the authors prove some estimates  for  $\lambda_F(\beta,\Omega)$ in terms of geometric quantities related to the domain $\Omega$, in particular, a  weighted anisotropic Cheeger inequality.

The aim of this paper is to prove, for a positive and continuous function $\beta$, two inequalities involving  $\lambda_F(\beta,\Omega)$ in terms  of   the following quantities
\begin{equation}
\label{vq}
\sigma_F(\beta,\Omega):= \inf_{\substack{v \in H^{1}(\Omega)\setminus \{0\}\\ \int_{\partial \Omega}\beta(x) vF(\nu)\,d \mathcal H^{n-1}=0}}\displaystyle \frac{\int_{\Omega}F^2(\nabla v)\,dx}{\int_{\Omega}v^2\,dx},
\end{equation}
and
\begin{equation}\label{qbis}
q_F(\beta,\Omega):= \inf_{\substack{h \in H^{1}(\Omega)\setminus \{0\}\\ \Delta_F h=0}}
\frac{\int_{\partial \Omega}\beta(x)h^2 F(\nu)d \mathcal H^{n-1}}{\int_{\Omega}h^2dx}.
\end{equation}
We observe that if $\beta(x)=\beta$ is a positive parameter, then
\begin{equation}\label{qbisc}
q_F(\beta,\Omega)=\beta q_F(\Omega),
\end{equation}
where
\begin{equation}
\label{qf}
 q_F(\Omega):= \inf_{\substack{h \in H^{1}(\Omega)\setminus \{0\}\\ \Delta_F h=0}}
\frac{\int_{\partial \Omega}h^2 F(\nu)d \mathcal H^{n-1}}{\int_{\Omega}h^2dx}
\end{equation}
while $\sigma_F( \beta,\Omega)$ does not depend on $\beta$ and then, in this case, we denote it by $\sigma_F(\Omega)$.
On the other hand, in the Euclidean case, when $\Omega$ has two axes of symmetry,  $\sigma_{\mathcal E}(\Omega)$ coincides with the first non-trivial Neumann eigenvalue of the Laplace operator $\mu(\Omega)$ (see for instance \cite{k}, \cite{ks1} and \cite{ks2}).  

Furthermore, under certain assumptions $q_{\mathcal E}(\Omega)$ coincides with the first nontrivial Steklov egenvalue $q$ related to the biharmonic Laplacian 
%
\begin{equation}
\label{P0}
\left\{
\begin{array}{ll}
\Delta^{2}v =0 &\text{in }\Omega \\
v= 0 & \text{on }\partial \Omega\\
\Delta v =q \dfrac{\partial v}{\partial \nu} & \text{on }\partial \Omega.
\end{array}
\right.
\end{equation}
This is shown in \cite{bfg} by means of a generalized Fichera's duality principle, provided $\Omega$ satisfies a uniform outer ball condition.
We recall that
 problem \eqref{P0} was first considered by Kuttler and Sigillito in \cite{ks1} and \cite{ks2} where among other things they studied the isoperimetric
properties related to the first eigenvalue. In the last years this kind of problems have been intensively studied in the literature, we refer the reader for instance to \cite{fgw}, \cite{s}, \cite{bgm}, \cite{bk} and the references therein for further studies. In particular  one can find some physical interpretation of the Steklov boundary conditions in \cite{bk} where the authors state also several Navier-Robin problems for the biharmonic operator.
Finally, when   $\beta$ is not a positive constant in \cite{hp}  the authors prove that  if $\partial \Omega \in C^{2}$ then $q_{\mathcal E}(\beta, \Omega)$ coincides with the first nontrivial eigenvalue of the following   ``weighted "  Steklov type problem
\begin{equation}
\label{var}
\left\{
\begin{array}{ll}
\Delta^{2}v =0 &\text{in }\Omega \\
v= 0 & \text{on }\partial \Omega\\
\Delta v = q_{\mathcal E} (\beta, \Omega )\dfrac{1}{\beta(x)} \dfrac{\partial v}{\partial \nu} & \text{on }\partial \Omega .
\end{array}
\right.
\end{equation}

Our main result is the following 
\begin{thm}
\label{main}
Let $\Omega \subset \R^n$, $n \ge 2$, be a bounded, connected, open set with Lipschitz boundary, then
\begin{equation}\label{Kutt_mi}
\frac{1}{\lambda_F(\beta,\Omega)}\leq\frac{1}{\sigma_F(\beta,\Omega)}+\frac{|\Omega|}{m}
\end{equation}
and
\begin{equation}\label{Kutt_40i}
\frac{1}{\lambda_F(\beta,\Omega)}\leq\frac{1}{\lambda_F(\Omega)}+\frac{1}{q_F(\beta,\Omega) },
\end{equation}
where $\beta\colon \partial\Omega \to ]0,+\infty[$ is a continuous  function, $m$ is given by \eqref{m} and $\lambda_F(\Omega)$ is the first Dirichlet eigenvalue of the Finsler Laplacian.
\end{thm}
In the Euclidean case, when $\beta(x)=\beta>0$ is  constant, the inequalities \eqref{Kutt_mi} and \eqref{Kutt_40i} were proved in  \cite{sp1} for  $\lambda_{\mathcal E}(\beta,\Omega)$, the first Robin eigenvalue of the Laplacian by using the P-function method (see also \cite{sp2}).

Successively in  \cite{k}, Kuttler  proves the same result with a simpler proof whose key ingredient is an algebraic inequality between geometric and arithmetic means.  Our result, in this order of idea, allows to extend  the results of  \cite{sp1} to the case when $\beta$ is not constant and to a larger class of elliptic operators. Our proof follows the idea contained in \cite{k}.\\
We  prove inequalities  \eqref{Kutt_mi} and   \eqref{Kutt_40i} in Section 2 and Section 3, respectively.

\section{Proof of the inequality \eqref{Kutt_mi}}
First of all, in order to verify that inequality \eqref{Kutt_mi} is well posed, we show that  $\sigma_F(\beta,\Omega)$ is positive.  To see this we first prove that $\sigma_F(\beta,\Omega)$ is a minimum.
 Let $v_k\in H^1(\Omega)\setminus\{0\}$ be a minimizing sequence such that $\int_{\partial \Omega}\beta(x) v_{k}F(\nu)\,d \mathcal H^{n-1}=0$, $\|v_{k}\|_{L^2(\Omega)}=1$
and
\begin{equation}
\label{lim} \lim_{k} \int_{\Omega}F^2(\nabla v_{k})\,dx=\sigma_F(\beta,\Omega).
\end{equation}
Then $v_k$  is bounded
in $H^1(\Omega)$ and there exists a subsequence, still denoted by $v_k$, such that $v_{k}$ converges in $L^{2}(\Omega)$  to a function  $v\in H^1(\Omega)$ with $\|v\|_{L^2(\Omega)}=1$. Furthermore, by the classical trace embedding Theorem,
 $v_k$ converges to $v$ also in $L^{2}(\partial\Omega)$ and  then $\int_{\partial \Omega}\beta(x) vF(\nu)\,d \mathcal H^{n-1}=0$. Taking $v$ as test function  in \eqref{vq} and using Fatou's Lemma, we  finally get \[
\sigma_F(\beta,\Omega) \le \int_{\Omega} F^2(\nabla v)\,dx\leq \liminf_{k} \int_{\Omega} F^2(\nabla v_k) \,dx= \sigma_F(\beta,\Omega).
\]
So  $v$ is a minimizer for $\sigma_F(\beta,\Omega)$ and  we can show that $\sigma_F(\beta,\Omega)>0$ arguing by contradiction. Indeed, if by absurd
$\sigma_F(\beta,\Omega)= 0$ then there exists $v\in H^1(\Omega)\setminus\{0\}$ such that $\int_{\partial \Omega} \beta(x) v \,F(\nu)\,d \haus=0$, $\|v\|_{L^{2}(\Omega)}=1$ and  $\int_\Omega F^2(\nabla u)\,dx=0$ and hence $v=C$ almost everywhere in $\overline \Omega$, with $C\in \mathbb R$, and $ C \neq 0$. Since $m>0$, this is in contradiction with
$$C\int_{\partial \Omega} \beta(x)  \,F(\nu)\,d \haus=0.$$
Now we  prove inequality \eqref{Kutt_mi}. 
For the reader's convenience from now on we will use the following notation 
\begin{equation}
\label{D}
E_F(w):=\displaystyle \int_{\Omega}F^2(\nabla w)\,dx,
\end{equation}
for any $w \in H^{1}(\Omega)$. Let $u$ be a positive eigenfunction corresponding to  $\lambda_F(\beta,\Omega)$
and 
\begin{equation}\label{c}
c=\dfrac{1}{m} \int_{\partial \Omega} \beta(x) u \,F(\nu)\,d \haus,
\end{equation}
with $m$ defined in \eqref{m}. By the Minkowski inequality and, recalling the definition of $\sigma_F(\beta,\Omega)$ in   \eqref{vq}, we have
\[
\sqrt{\int_{\Omega}u^2dx}\leq \sqrt{\int_{\Omega}(u-c)^2}+\sqrt{c^2 |\Omega|}\leq \sqrt{\frac{E_F(u)}{\sigma_F(\beta,\Omega)}}+\sqrt{c^2 |\Omega|}.
\]
Squaring and using the arithmetic-geometric mean inequality, we have
\begin{equation}\label{A}
\begin{split}
\int_{\Omega}u^2dx &\leq \frac{E_F(u)}{\sigma_F(\beta,\Omega)}+c^2 |\Omega|+2\sqrt{\frac{E_F(u)c^2|\Omega|}{\sigma_F(\beta,\Omega)}}\\
&\leq \frac{E_F(u)}{\sigma_F(\beta,\Omega)}+c^2 |\Omega|+\frac{E_F(u)|\Omega|}{m}+\frac{c^2 m}{\sigma_F(\beta,\Omega)}\\
&= E_F(u)\left(\frac{1}{\sigma_F(\beta,\Omega)}+\frac{|\Omega|}{m}\right)+c^2m\left(\frac{1}{\sigma_F(\beta,\Omega)}+\frac{|\Omega|}{m}\right)\\
&=\left(\frac{1}{\sigma_F(\beta,\Omega)}+\frac{|\Omega|}{m}\right)(E_F(u)+c^2 m).
\end{split}
\end{equation}
By \eqref{c}, H\"{o}lder inequality and \eqref{lambda}, we see that \eqref{A} implies 
\begin{equation}
\begin{split}
\int_{\Omega}u^2dx &\leq \left(\frac{1}{\sigma_F(\beta,\Omega)}+\frac{|\Omega|}{m}\right)\left(E_F(u)+
\frac{\left(\int_{\partial \Omega}\beta(x) F(\nu)d \mathcal H^{n-1}\right)\left(\int_{\partial \Omega}\beta(x)u^2 F(\nu)d \mathcal H^{n-1}\right)}{m}\right)\\
&=\left(\frac{1}{\sigma_F(\beta,\Omega)}+\frac{|\Omega|}{m}\right)\left(E_F(u)+\int_{\partial \Omega}\beta(x)u^2 F(\nu)d \mathcal H^{n-1}\right)\\
&=\left(\frac{1}{\sigma_F(\beta,\Omega)}+\frac{|\Omega|}{m}\right)\left(\lambda_F(\beta,\Omega)\int_{ \Omega}u^2 dx\right)
\end{split}
\end{equation}
which gives \eqref{Kutt_mi}.
\begin{rem}
Let $\Omega$ be an open set of $\R^{n}$ with Lipschitz boundary. We denote by $P_{F}(\Omega)$ the so-called anisotropic perimeter defined as follows (see for instance \cite{ab})
\[
P_F(\Omega)=\displaystyle \int_{\de \Omega}F(\nu) \, d \mathcal H^{n-1},
\]
where $\nu$ denotes the unit outer normal to $\de \Omega$. We stress that when $\beta(x)=\beta$ is a positive constant, the inequality \eqref{Kutt_mi} gives the following  asymptotic behavior of $\lambda_F(\beta,\Omega)$, when $\beta$ goes to zero:
\begin{equation}
\label{asi}
\lim_{\beta \rightarrow0} \frac{\lambda_F(\beta,\Omega)}{\beta}=\frac{P_{F}(\Omega)}{|\Omega|}.
\end{equation}
Indeed if $\beta$ is a positive constant then  $m=\beta P_F(\Omega)$ and 
 we have
\begin{equation}
\label{s1}
\frac{P_{F}(\Omega)}{|\Omega|} \ge \frac{\lambda_F(\beta,\Omega)}{\beta}\ge \dfrac{P_{F}(\Omega)\sigma_F(\Omega)}{P_{F}(\Omega) \beta + |\Omega|\sigma_F(\Omega)},
\end{equation}
where the first inequality follows by  using a constant as test function in \eqref{lambda} and the second by using \eqref{Kutt_mi}. Taking in  \eqref{s1} the limit for $\beta $ which goes to zero  one get \eqref{asi}.
\end{rem}
\begin{rem}
Let $\mu_F(\Omega)$ be the first non-trivial Neumann eigenvalue of  the Finsler Laplacian (see for instance \cite{dpgp}), it holds
\begin{equation}
\label{confr}
\sigma_F(\beta,\Omega) \le \mu_F(\Omega).
\end{equation}
Indeed if $u$ is an eigenfunction corresponding to $\mu_F(\Omega)$, then $\int_{\Omega}u\,dx=0$ and
\[
\mu_F(\Omega)=\displaystyle \frac{\int_{\Omega}F^2(\nabla u)\,dx}{\int_{\Omega}u^2\,dx}.
\]
Inequality \eqref{confr} follows by taking   as test  in \eqref{vq} the function $v(x)=u(x)-c$, where $c$ is as in \eqref{c}.
\end{rem}
\section{Proof of the inequality \eqref{Kutt_40i}}
First of all, we observe that the  trace embedding Theorem ensures that   $q_F(\beta,\Omega)$ is positive and then inequality \eqref{Kutt_40i} is well posed.\\
Let $u$ be a positive eigenfunction corresponding to  $\lambda_F(\beta,\Omega)$ and let us consider the functions $v$ and $h$ which solve the following problems respectively
\begin{equation}
\label{v}
\left\{
\begin{array}{ll}
\Delta_Fv=\Delta_Fu & \text{in }\Omega \\
v=0 & \text{on }\partial \Omega ,
\end{array}
\right.
\end{equation}
and
\begin{equation}
\label{h}
\left\{
\begin{array}{ll}
\Delta_Fh=0 & \text{in }\Omega \\
h=u & \text{on }\partial \Omega.
\end{array}
\right.
\end{equation}
The maximum principle assures that $u \le v+h$. Moreover, by using the same notation of Section 2,   it holds
\begin{equation}
\label{grad}
E_F(u)\ge E_F(v).
\end{equation}
Indeed, the convexity of $F^{2}$ and the homogeneity of  $F$ imply
\begin{align*}
\int_{\Omega} F^2(\nabla v) dx &\geq -\int_{\Omega}F^2(\nabla u) dx+2\int_{\Omega}F(\nabla v)F_{\xi}(\nabla v)\cdot \nabla v dx
\\& = -\int_{\Omega}F^2(\nabla u) dx +2\int_{\Omega}F^2(\nabla v) dx,
\end{align*}
where last equality follows being $v$ the solution of \eqref{v}. \\ By the Minkowski  inequality, by  the definition of $q_F(\beta,\Omega)$ given in \eqref{qbis} and recalling the following variational characterization of $\lambda_F(\Omega)$ (see for instance \cite{bfk})
\begin{equation}
\lambda_{F}(\Omega)=\min_{u \in H_{0}^{1}(\Omega)\setminus \{0\}}\displaystyle \frac{E_{F}(u)}{\displaystyle\int_{\Omega}u^{2}\,dx},
\end{equation}
 we get
\begin{equation*}
\begin{split}
\left(\int_{\Omega}u^2dx\right)^{\frac{1}{2}} &\leq \left(\int_{\Omega}v^2dx\right)^{\frac{1}{2}}+\left(\int_{\Omega}h^2dx\right)^{\frac{1}{2}}\\
&\leq\left(\frac{E_F(v)}{\lambda_F(\Omega)}\right)^{\frac{1}{2}}+\left(\frac{\int_{\partial \Omega}\beta(x)h^2 F(\nu)d \mathcal H^{n-1}}{q_F(\beta,\Omega)}\right)^{\frac{1}{2}}\\
&\leq\left(\frac{E_F(u)}{\lambda_F(\Omega)}\right)^{\frac{1}{2}}+\left(\frac{\int_{\partial \Omega}\beta(x)u^2 F(\nu)d \mathcal H^{n-1}}{q_F(\beta,\Omega)}\right)^{\frac{1}{2}}.
\end{split}
\end{equation*}
Squaring and using the arithmetic-geometric mean inequality, we have
\begin{equation*}
\begin{split}
\int_{\Omega}u^2dx &\leq \frac{E_F(u)}{\lambda_F(\Omega)}+\frac{\int_{\partial \Omega}\beta(x)u^2 F(\nu)d \mathcal H^{n-1}}{q_F(\beta,\Omega)}+2\left(\frac{E_F(u)\int_{\partial \Omega}\beta(x)u^2 F(\nu)d \mathcal H^{n-1}}{\lambda_F(\Omega)q_F(\beta,\Omega)}\right)^{\frac{1}{2}}\\
&\leq\left(\frac{1}{\lambda_F(\Omega)}+\frac{1}{q_F(\beta, \Omega)}\right)E_F(u)+\int_{\partial \Omega}\beta(x)u^2 F(\nu)d \mathcal H^{n-1}\left(\frac{1}{q_F(\beta,\Omega)}+\frac{1}{\lambda_F(\Omega)}\right)\\
&\leq\left(\frac{1}{\lambda_F(\Omega)}+\frac{1}{q_F(\beta, \Omega)}\right)\left(E_F(u)+\int_{\partial \Omega}\beta(x)u^2 F(\nu)d \mathcal H^{n-1}\right),
\end{split}
\end{equation*}
which gives \eqref{Kutt_40i}.
\begin{rem}
We stress that inequality \eqref{Kutt_40i} can also be written as follows
\begin{equation}\label{cont1}
0\le\frac{1}{\lambda_F(\beta,\Omega)}-\frac{1}{\lambda_F(\Omega)}\le\frac{1}{q_F(\beta,\Omega) },
\end{equation}
this inequality gives an upper bound of the distance between the first Dirichlet and Robin eigenvalue of the Finsler Laplacian in terms of $q_F(\beta,\Omega) $. In particular, if $\beta(x)=\beta$ is constant, \eqref{cont1} can be rewritten by \eqref{qbisc} as
\begin{equation}\label{cont2}
0\le\frac{1}{\lambda_F(\beta,\Omega)}-\frac{1}{\lambda_F(\Omega)}\le\frac{1}{\beta q_{F}(\Omega) },
\end{equation}
with $q_{F}(\Omega)$ defined by \eqref{qf}. Inequality \eqref{cont2} implies that if $\beta \to +\infty $ then $\lambda_F(\beta,\Omega) \to \lambda_F(\Omega)$ as well known.
\end{rem}
\begin{rem}
We observe that \eqref{Kutt_40i} gives a geometric inequality involving $q_F(\beta,\Omega)$. Indeed,  if we consider $h\equiv1$ as test function in \eqref{qbis} we have
\begin{equation}\label{cont}
q_F(\beta,\Omega)\le \frac{m}{|\Omega|},
\end{equation}
that reads in the constant case by \eqref{qbisc} as
\[
q_F(\Omega)\le \frac{P_F(\Omega)}{|\Omega|},\]
where $q_{F}(\Omega)$ is defined in \eqref{qf}.
\end{rem}
\section*{Acnowledgements}

This work has been partially supported by the MiUR-PRIN 2017 grant \lq\lq Qualitative and quantitative aspects of nonlinear PDEs\rq\rq, by GNAMPA of INdAM, by FRA 2020 \lq\lq Optimization problems in Geometric-functional inequalities and nonlinear  PDEs\rq\rq (OPtImIzE), by project Vain-Hopes within the program \\VALERE:VAnviteLli pEr la RicErca.


\begin{thebibliography}{20}

\bibitem{aflt}
A. Alvino, V. Ferone, P.-L. Lions, G. Trombetti,
\textit{Convex symmetrization and applications},
Ann. Inst. H. Poincar\'{e} Anal. Non Lin\'{e}aire, 14 (1997) 275--293.

\bibitem{ab}
M.~Amar and G.~Bellettini.
\newblock \textit{A notion of total variation depending on a metric with discontinuous coefficients}.
\newblock {Ann. Inst. H. Poincar\'{e} Anal. Non Lin\'{e}aire},
  11 (1994) 91--133.
\bibitem{bfk}  M. Belloni, V. Ferone, B. Kawohl, \emph{Isoperimetric inequalities, Wulff shape and related questions for
strongly nonlinear elliptic operators}, Z. Angew. Math. Phys. ZAMP 54 (2003) 771-783.

\bibitem{fk} V. Ferone, B. Kawohl, \emph{Remarks on a Finsler-Laplacian}, Proc. Amer. Math. Soc. 137 (1) (2009) 247-253.


\bibitem{bgm}
E. Berchio, F. Gazzola, E. Mitidieri, \textit{ Positivity preserving property for a class of biharmonic elliptic
problems}, J. Differ. Equa. 229 (2006) 1-23.


\bibitem{bfg}
D. Bucur, A. Ferrero, F. Gazzola,
\textit{On the first eigenvalue of a fourth order Steklov problem}, Calc. Var., 35 (2009) 103-131.

\bibitem{bk}
D. Buoso, J.B. Kennedy,
\textit{The Bilaplacian with Robin boundary conditions}, preprint 2021, arXiv:2105.11249 .

 \bibitem{dpdfg} F. Della Pietra, G. di Blasio, N. Gavitone, \textit{Sharp estimates on the first Dirichlet eigenvalue of nonlinear elliptic operators via maximum principle}, Adv. Nonlinear Anal. 9 (2020) 278--291.
\bibitem{pota} F. Della Pietra, N. Gavitone, \textit{Faber-Krahn Inequality for Anisotropic Eigenvalue
Problems with Robin Boundary Conditions}, Potential Anal 41 (2014) 1147--1166.

\bibitem{dpgp} F. Della Pietra, N. Gavitone, G. Piscitelli, \textit{A sharp weighted anisotropic Poincar\'e inequality for convex domains}, Comptes Rendus Mathematique 355 ( 2017) 748--752.

\bibitem{dgn} F. Della Pietra, N. Gavitone \emph{Anisotropic elliptic problems involving Hardy-type potentials}, J. Math. Anal. Appl. 397 (2013) 800-813.

\bibitem{fgw} A. Ferrero, F. Gazzola, T. Weth, \textit{On a fourth order Steklov eigenvalue problem}, Analysis 25 (2005) 315-332.

\bibitem{f}G. Fichera, \emph{Su un principio di dualit\`a per talune formule di maggiorazione relative alle equazioni differenziali}, Atti Accad. Naz. Lincei 19 (1955) 411-418.
\bibitem{gt} N. Gavitone, L. Trani, \emph{On the first Robin eigenvalue of a class of anisotropic operators}, Milan J. Math. 86 (2018) 201-223.
\bibitem{hp}J. Hersch, L.E. Payne,\emph{One-Dimensional Auxiliary Problems and a priori Bounds}, Abhandlungen Math. Seminar Univ. Hamburg, Band 36 (1971) 57-65.
\bibitem{k} J. R. Kuttler, \emph{A Note on a Paper of Sperb}. ZAMP 24 (1973) 431-434.
\bibitem{ks1} J.R. Kuttler, V. G. Sigillito, \emph{Inequalities for Membrane and Stekloff Eigenvalues}, J. Math. Anal. Appl. 23 (1968) 148-160.
\bibitem{ks2} J.R. Kuttler, V. G. Sigillito, \emph{An Inequality for a Stekloff Eigenvalue by the Method of Defect}, Proc. Amer. Math. Soc. 20 (1969) 357--360.
\bibitem{p} L.E. Payne,\emph{Isoperimetric Inequalities and their Applications}, SIAM Rev. 9 (1967) 453--488.
\bibitem{s} S. Raulot, A. Savo, \textit{Sharp Bounds for the First Eigenvalue of a Fourth-Order Steklov Problem},
J. Geom. Anal., 25 (2015) 1602--1619.
\bibitem{sp1} R. P. Sperb,\textit{Untere und obere Schranken fur den tiefsten Eigenwert der elastisch gestutzten Membran},
Z. Angew. Math. Phys. 23  (1972) 231--244.
\bibitem{sp2} R. P. Sperb, \textit{Maximum principles and applications}, Academic Press, 1981.



\end{thebibliography}
\end{document}